\documentclass[12pt]{article}
\usepackage{amscd}
\usepackage{mathrsfs}
\usepackage{amsfonts}
\usepackage{graphicx}
\usepackage{amsmath,amscd,stmaryrd,amssymb,array, amsfonts,amsmath,amssymb}
\usepackage{enumitem}
\setenumerate[1]{itemsep=0pt,partopsep=0pt,parsep=\parskip,topsep=5pt}
\setitemize[1]{itemsep=0pt,partopsep=0pt,parsep=\parskip,topsep=0pt}
\setdescription{itemsep=0pt,partopsep=0pt,parsep=\parskip,topsep=5pt}

\numberwithin{figure}{section}
 \numberwithin{equation}{section}

\newtheorem{theorem}{Theorem}[section]
\newtheorem{proposition}[theorem]{Proposition}
\newtheorem{definition}[theorem]{Definition}

\newtheorem{lemma}[theorem]{Lemma}
\newtheorem{remark}[theorem]{Remark}

\newcommand{\cA}{{\mathcal A}}

\newcommand{\cM}{{\mathcal M}}

\newcommand{\cT}{{\mathcal T}}

\newcommand{\sX}{{\mathscr X}}

\newcommand{\mB}{\mb{B}}

\def\be{\begin{equation}}
\def\ee{\end{equation}}
\def\bes{\begin{equation*}}
\def\ees{\end{equation*}}
\def\bsp{\begin{split}}
\def\esp{\end{split}}

\def\ba{\begin{array}}
\def\ea{\end{array}}
\def\benu{\begin{enumerate}}
\def\eenu{\end{enumerate}}
\def\bt{\begin{theorem}}
\def\et{\end{theorem}}
\def\bp{\begin{proposition}}
\def\ep{\end{proposition}}
\def\bl{\begin{lemma}}
\def\el{\end{lemma}}
\def\br{\begin{remark}}
\def\er{\end{remark}}
\def\bd{\begin{definition}}
\def\ed{\end{definition}}

\def\b{\beta}
\def\De{\Delta}
\def\de{\delta}

\def\nab{\nabla}
\def\lam{\lambda}
\def\Lam{\Lambda}

\def\ve{\varepsilon}
\def\sig{\sigma}

\def\gam{\gamma}

\def\a{\alpha}
\def\W{\Omega}

\def\.{\cdot}

\def\R{\mathbb{R}}

\def\A{\forall}
\def\ol{\overline}
\def\ul{\underline}
\def\Cap{\bigcap}\def\Cup{\bigcup}

\def\ra{\rightarrow}

\def\~{\tilde}
\def\8{\infty}
\def\X{\times}
\def\({\left(}
\def\){\right)}

\def\mb{\mbox}

\def\Hs{\hspace{1cm}}\def\hs{\hspace{0.5cm}}
\def\Vs{\vskip8pt}\def\vs{\vskip4pt}

\def\({\left(}\def\){\right)}

\begin{document}

\begin{center}{\large BIFURCATION FROM INFINITY AND MULTIPLICITY\\ OF THE LANDESMAN-LAZER TYPE PROBLEM\\ [1ex]OF ELLIPTIC EQUATIONS}
\end{center}
\Vs
\begin{center}Xuewei Ju
\end{center}

\begin{center}
{\em Department of Mathematics, Civil Aviation University of China\\
     Tianjin 300300,  China }

\begin{center}Desheng Li$^*$\footnote{$*$ Corresponding author: Desheng Li}\, and\, Youbin Xiong
\end{center}

{\em  School of Mathematics,  Tianjin University\\
  Tianjin 300072, China }

  \Vs
  (xwju@cauc.edu.cn; lidsmath@tju.edu.cn)
\end{center}

\Vs

{\em Abstract}\,\,  This paper deals  with the Landesman-Lazer type problem of elliptic equations associated with homogeneous Dirichlet boundary conditions. By using some dynamical arguments  we derive some new results on bifurcation from infinity and multiplicity of the problems.


\section{Introduction}
In this paper we consider the Landesman-Lazer type problem for the boundary value problem:
\be\label{eq1}\left\{\ba{ll}
-\Delta u=\lam u+f(x,u), \hs\hs x\in
  \Omega;\\[1ex]
u(x)=0, \Hs\hs\, \Hs\Hs x\in\partial\Omega
\ea\right. \ee as $\lam$ varies near resonance, where  $\W\subset\R^N$ is a bounded domain, and $f\in C^1(\ol\Omega\X\R)$  satisfies the Landesman-Lazer type condition:
\benu
 \item[(LC)] $f$ is bounded; furthermore,
\be\label{LC1}   \liminf\limits_{t\rightarrow+\infty}f(x,t)\geq\ol{f}>0,\Hs
  \limsup\limits_{t\rightarrow-\infty}f(x,t)\leq-\underline{f}<0\ee
  uniformly for $x\in\ol\W$ (where $\ol f$ and $\ul f$ are independent of $x$).
  \eenu
Such problems can be seen as nonlinear perturbations of the corresponding linear ones, and has aroused much interest in the past decades; see  \cite{CL,CM,FGP,LW,PM,PP,ST,H.N,BC,BEG,K,CC,U,MS1,MS,SW,CW} and references therein.

If $\lam$ is not an eigenvalue of the operator $A=-\De$ (associated with the homogeneous Dirichlet boundary condition), it can be easily shown that the solution set of \eqref{eq1} is bounded. This basic fact in turn allows us to show that the problem has at least a solution by using different means, in particular by means of fixed point theory and topological degree. 
Here we are interested in the multiplicity of solutions of \eqref{eq1} as $\lam$ varies near an eigenvalue $\mu_k$ of $A$. The motivation comes from the work of  Mawhin and Schmitt \cite{MS1,MS}, Schmitt and Wang \cite{SW} and Chang and Wang \cite{CW}, etc.

  In \cite{MS} the authors proved under appropriate Landesman-Lazer type conditions that if $\mu_k$ is of odd multiplicity, then the problem has  at least two distinct solutions for $\lam$ on one side of $\mu_k$ but close to $\mu_k$, and at least one solution for $\lam$ on the other side.  Later the restriction on the multiplicity of $\mu_k$ in this result was removed by Schmitt and Wang
in a general framework on bifurcation of potential operators in \cite{SW}. A pure dynamical argument for  Schmitt and Wang's result on \eqref{eq1} can also be found in \cite{LLZ}.

For the first eigenvalue $\mu_1$, it was shown in \cite{MS1} (see \cite[Theorems 4 and 5]{MS1}) and \cite{CW} (see \cite[Section 3]{CW}) that the problem has at least three distinct solutions for $\lam$ in a one-sided neighborhood of $\mu_1$, two of which going to infinity and one remaining bounded as $\lam\ra\mu_1$. Our main purpose in this present work is to extend this elegant result to any eigenvalue  $\mu_k$ of $A$.
Specifically, let $$
\b_k=\min\{\mu_k-\mu_{k-1},\mu_{k+1}-\mu_k\},\Hs k=1,2,\cdots
$$
(here we assign $\b_1=\mu_2-\mu_1$), and assume that
\be\label{SMC}M\mu_1^{-1/2}L_f\int_0^\8\big(2+\tau^{-\frac{1}{2}}\big)e^{-\frac{1}{4}\b_k\tau}d\tau<1,\ee
 where   $M\geq1$ is a constant depending only upon the operator $A$, and $L_f$ is the Lipschitz constant of $f$. We will show under the above smallness requirement on  $L_f$ that  there exists $0<\theta \leq \b_k/4$ such that \eqref{eq1} has at least three distinct solutions $u_{\lam}^{i}$ ($i=0,1,2$) for each $\lam\in [\mu_k-\theta ,\mu_k)$, and $u_{\lam}^{1}$ and $u_{\lam}^{2}$ go to $\8$ as $\lam\ra \mu_k$ whereas $u_\lam^0$ remains bounded on $[\mu_k-\theta ,\mu_k)$.
A ``dual'' version of the result also holds true if we replace \eqref{LC1} by
\be\label{LC2} \limsup\limits_{t\rightarrow +\infty}f(x,t)\leq-\ol{f}<0,\Hs   \liminf\limits_{t\rightarrow -\infty}f(x,t)\geq\ul{f}>0.\ee

It is worth noticing that for a given globally Lipschitz continuous  function $f$,  since $\b_k\ra+\8$ as $k\ra\8$ (hence the integral in \eqref{SMC} goes to $0$), the condition \eqref{SMC} is automatically fulfilled provided $k$ is sufficiently large.

Our method here is as follows. Instead of transforming \eqref{eq1} into an operator equation and applying the topological degree or other means such as variational methods, we view the problem as the stationary one of  the parabolic equation \be\label{eq02}\left\{\ba{ll}
  u_t-\Delta u=\lam u+f(x,u), \hs\,\, x\in
  \Omega;\\[1ex]
  \,u(x)=0, \Hs\Hs \Hs \Hs x\in\partial\Omega.
\ea\right. \ee Then we give  an existence result on some global invariant manifolds $\cM_\lam^c$ for the semiflow $\Phi_\lam$ generated by \eqref{eq02} for $\lam$ near each eigenvalue $\mu_k$. Such a  manifold $\cM_\lam^c$ contains all the invariant sets of the system. This allows us to reduce the system on $\cM_\lam^c$ and prove, using  the shape theory of attractors \cite{Kap}, that there exists  $0<\theta \leq \b_k/4$ such that the system  bifurcates from infinity a compact isolated invariant set $K_\lam^\8$ which takes the shape of a sphere $\mathbb{S}^{m-1}$ for $\lam\in\Lam_k^-=[\mu_k-\theta,\mu_k)$, where $m$ is the algebraic multiplicity  of $\mu_k$. Since $\Phi_\lam$ is a gradient system, it can be shown that $K_\lam^\8$ necessarily contains two distinct equilibria of $\Phi_\lam$. These equilibria are precisely solutions of \eqref{eq1}. Thus we conclude  that the Landesman-Lazer type problem \eqref{eq1} bifurcates from infinity two distinct solutions as $\lam$ varies in $\Lam_k^-$.
 Combining this result with some known ones in \cite[Theorem 5.2]{LLZ}, we immediately complete the proof of  our main results promised above.

Let us mention that our approach is of a pure dynamical nature and is different from those in the literature. It  allows  us to obtain a more clear picture on the dynamic bifurcation from infinity of the parabolic problem \eqref{eq02} near resonance, which, from the point of view of dynamical systems theory, is naturally of independent interest.

This paper is organized as follows. In Section 2 we present an existence result on  global invariant   manifolds of \eqref{eq02} in an abstract framework of evolution equations in Banach spaces. In Section 3 we give a more precise description on the dynamic bifurcation from infinity of \eqref{eq02} and prove our main results.
\section{Existence of invariant manifolds for nonlinear evolution equations} Let $X$ be a Banach space with norm $\|\.\|$, and $A$ be a sectorial operator on $X$ with compact resolvent.  Consider the semilinear equation \be\label{e:2.0}x_t+Ax=\lam x+f(x)\ee
in $X$. In this section we present an existence result on global invariant manifolds for the equation when $\lam$ varies near the real part  $\lam_0=\hbox{Re}\,\mu_0$ of an eigenvalue  $\mu_0$ of the operator $A$.

\subsection{Mathematical setting}
Denote $\sig(A)$ the spectrum of $A$ and write
$\mb{Re}\,\sig(A):=\{\mb{Re}\,\mu:\,\,\mu\in \sig(A)\}.$
  Pick a number $a>0$ such that $$\hbox{Re}\,\sig(A+aI)>0.$$
Let  $\Lambda=A+aI$. For each $\a\geq 0$, define
$X^\a=D(\Lambda^\a)$. $X^\a$ is equipped {with \,\,the\,\,norm} $\|\.\|_\a$ defined by  $$\|x\|_\a=\|\Lambda^\a x\|,\hs x\in X^\a.$$
It is well known  that the definition of $X^\a$ is independent of the choice of  $a$.

Let $\lam_0\in\mb{Re}\,\sig(A)$. Since $A$ has compact resolvent, $\lam_0$ is isolated in $\mb{Re}\,\sig(A)$. Hence   $\sig(A)$ has a spectral decomposition $\sig(A)=\sig_1\cup \sig_2\cup\sig_3$ with $\sig_2=\{\mu\in\sig(A):\,\,\hbox{Re}\,\mu=\lam_0\}$ and
$$\sig_1=\{\mu\in\sig(A):\,\,\hbox{Re}\,\mu<\lam_0\},\hs  \sig_3=\{\mu\in \sig(A):\,\,\hbox{Re}\,\mu>\lam_0\}.$$
Let
$$
\gam_1=\max\{\mb{Re}\,\lam:\,\,\lam\in\sig_1\},\hs \gam_3=\min\{\mb{Re}\,\lam:\,\,\lam\in\sig_3\}.
$$
Then $\gam_1<\lam_0<\gam_3$.

The space $X$ has a corresponding direct sum  decomposition $X=X_1\oplus X_2\oplus X_3$ with $X_1$
 and $X_2$ being  finite dimensional.
Denote $\Pi_i:X\ra X_i$   the projection from $X$ to $X_i$ ($i=1,2,3$).


Let $\b=\min\{\lam_0-\gamma_1,\,\gamma_2-\lam_0\},$ and write
$$B=B(\lam):=A-\lam I,\hs B_i=B|_{X_i}.$$ We infer from  Henry \cite[Theorems 1.5.3 and 1.5.4]{D.H} that  there exists $M\geq1$ (depending only upon $A$) such that for $\a\in[0,1]$
\be\label{b1}\|e^{-B_1 t}\|\leq Me^{\frac{3}{4}\b t},\hs \|\Lambda^\a e^{-B_1 t}\|\leq Me^{\frac{3}{4}\b t},\Hs t\leq0,\ee
\be\label{b3}\| e^{-B_2 t}\|\leq Me^{\frac{\b}{4}|t|},\hs \|\Lambda^\a e^{-B_2 t}\|\leq Me^{\frac{\b}{4}|t|},\Hs t\in\R.\ee
\be\label{b2}\|\Lambda^\a e^{-B_3 t}\|\leq Mt^{-\a}e^{-\frac{3}{4}\b t},\hs \|\Lambda^\a e^{-B_3 t}\Pi_3\Lambda^{-\a}\|\leq Me^{-\frac{3}{4}\b t},\Hs t>0,\ee
(The latter  estimates in \eqref{b1} and \eqref{b3} are due to the finite dimensionality of the spaces $X_1$ and $X_2$.)

Given $\mu\geq 0$, define a Banach space $\sX_\mu$ as
$$\mathscr{X}_\mu=\left\{x\in C(\R;X^\a): \,\,\sup_{t\in\R}e^{-\mu |t|}\|x(t)\|_\a<\8\right\},$$  which is equipped with the norm $\|\.\|_{\mathscr{X}_\mu}$: $$\|x\|_{\mathscr{X}_\mu}=\sup_{t\in\R}e^{-\mu |t|}\|x(t)\|_\a,\Hs \A\,x\in \sX_\mu.$$

The equation \eqref{e:2.0} can be rewritten as  \be\label{e1}x_t+B x=f(x).\ee For our purposes here, from now on we always assume
\benu
\item[(F1)] $f\in C(X^\a,X)$ and is  globally Lipschitz for some $\a\in[0,1)$.
\eenu
It is easy to see that this condition also  implies that there exists $C>0$ such that
\be\label{eC}
\|f(x)\|\leq C(\|x\|_\a+1),\Hs\A\,x\in X^\a.
\ee

  Under the assumption (F1)  the Cauchy problem of \eqref{e1} is well-posed in $X^\a$. Specifically, for each $x_0\in X^\a$  the equation \eqref{e1} has a unique strong solution $x(t)=\phi_\lam(t;x_0)$ with initial value $x(0)=x_0$ which globally exists on $\R^+$;  see, e.g., \cite[Corollary 3.3.5]{D.H}.
   Set $\Phi_\lam(t)x_0=\phi_\lam(t;x_0)$. Then $\Phi_\lam$ is a (global) semiflow on $X^\a$.

\subsection{A basic lemma }
The following lemma  will play a fundamental role in the construction of invariant manifolds.
\bl\label{le4.2} Let $\b/4<\mu<3\b/4$. A function $x\in\mathscr{X}_\mu$ is a solution of \eqref{e1} on $\R$ if and only if it solves the following integral equation
\begin{equation}\label{e4}
\begin{split}
x(t)&=e^{-B_2 t}\Pi_2 x(0)+\int_{0}^te^{-B_2(t-\tau)}\Pi_2f(x(\tau))d\tau\\
&\quad+\int_{-\8}^{t}e^{-B_3(t-\tau)}\Pi_3 f(x(\tau))d\tau\\
&\quad-\int_{t}^\8e^{-B_1(t-\tau)}\Pi_1 f(x(\tau))d\tau.
\end{split}
\end{equation}
\el
\br In case $\sig_1=\emptyset$ the integral equation \eqref{e4} reduces to \begin{equation*}
\begin{split}
x(t)&=e^{-B_2 t}\Pi_2 x(0)+\int_{0}^te^{-B_2(t-\tau)}\Pi_2 f(x(\tau))d\tau\\
&\quad+\int_{-\8}^{t}e^{-B_3(t-\tau)}\Pi_3 f(x(\tau))d\tau.\end{split}
\end{equation*}\er
{\bf Proof of Lemma \ref{le4.2}.}  Let $x\in \mathscr{X}_\mu$ be a solution of $(\ref{e1})$ on $\R$. Write  $x=x_1+x_2+x_3$, where $x_i=\Pi_i x$ ($i=1,2,3$). Then for any $t,t_0\in\R$,
\be\label{e6b} x_i(t)=e^{-B_i(t-t_0)}x_i(t_0)+\int_{t_0}^te^{-B_i(t-\tau)}\Pi_i f(x(\tau))d\tau,\hs i=1,2,3.\ee

For $i=1$, if $t\leq t_0$ then by \eqref{b3} we see that
\begin{equation}\label{e8}
\begin{split}
\|e^{-B_1(t-t_0)}x_1(t_0)\|_\a&\leq Me^{\frac{3}{4}\b(t-t_0)}\|x(t_0)\|_\a\\
&=Me^{\frac{3}{4}\b t}e^{-(\frac{3}{4}\b-\mu)t_0}e^{-\mu t_0} \|x(t_0)\|_\a\\
& \leq Me^{\frac{3}{4}\b t}e^{-(\frac{3}{4}\b-\mu)t_0}\|x\|_{\mathscr{X}_\mu}\ra0, \hs\hbox{as }t_0\ra+\8.\end{split}
\end{equation}
Thus setting $t_0\ra+\8$ in \eqref{e6b} we obtain that  $$x_1(t)=-\int^\8_te^{-B_1(t-\tau)}\Pi_1 f(x(\tau))d\tau.$$

For $i=3$, if $t_0\leq 0$ then  by \eqref{b2} we have
\begin{equation*}
\begin{split}
\|e^{-B_3(t-t_0)}x_3(t_0)\|_\a&\leq Me^{-\frac{3}{4}\b (t-t_0)} \|x(t_0)\|_\a\\
&=Me^{-\frac{3}{4}\b t}e^{(\frac{3}{4}\b-\mu)t_0}\(e^{\mu t_0} \|x(t_0)\|_\a\)\\
&\leq Me^{-\frac{3}{4}\b t}e^{(\frac{3}{4}\b-\mu)t_0}\|x\|_{\mathscr{X}_\mu}\ra0,\hs \mb{as }t_0\ra-\8.
\end{split}
\end{equation*}
Hence by \eqref{e6b} we deduce that
$$x_3(t)=\int^{t}_{-\8}e^{-B_3(t-\tau)}\Pi_3 f(x(\tau))d\tau.$$

For $i=2$, taking $t_0=0$ in \eqref{b2} it yields
$$
x_2(t)=e^{-B_2t}x_2(0)+\int_{0}^te^{-B_2(t-\tau)}\Pi_2 f(x(\tau))d\tau
$$
Combing the above results together one immediately concludes the validity of the equation \eqref{e4}.

Conversely, for each  $x\in\mathscr{X}_\mu$ satisfying $(\ref{e4})$ one can easily verify that $x$  is a solution of $(\ref{e1})$ on $\R$.  $\Box$

\subsection{Existence of global invariant  manifolds}
We are now ready to state and prove our main result in this section.

Denote  $$X^\a_i:=X_i\cap X^\a,\Hs  i=1,2,3,$$ and let $X^\a_{ij}=X^\a_i\oplus X^\a_j$ ($i\ne j$).
\bt\label{th1} Suppose that the Lipschitz constant $L_f$ of $f$ satisfies
 \be\label{eq1.6}ML_f\int_0^\8\big(2+\tau^{-\a}\big)e^{-\frac{\b}{4}\tau}d\tau <1.\ee
 Then for each $\lam\in (\lam_0-\b/4,\lam_0+\b/4):=J$,
the semiflow $\Phi_\lam$  defined by \eqref{e:2.0} has a global invariant manifold $\cM^c_\lam$ given by $$\cM^c_\lam=\{y+\xi_\lam(y):y\in X^\a_2\},$$ where $\xi_\lam:X^\a_2\ra X^\a_{13}$ is Lipschitz continuous uniformly on $\lam\in J$.
\et

\Vs\noindent{\bf Proof of Theorem \ref{th1}.}  Instead of the original equation we consider the modified one in \eqref{e1}. For each $\lam\in J$ and $y\in X_2^\a$, one can use the righthand side of equation $(\ref{e4})$ to define a contraction mapping $\cT:=\cT_{\lam,y}$ on
 $\mathscr{X}_{\b/2}$ as follows:
\begin{equation*}
\begin{split}
\cT x(t)&=e^{-B_2 t}y+\int_{0}^te^{-B_2(t-\tau)}\Pi_2 f(x(\tau))d\tau\\
&\quad+\int_{-\8}^{t}e^{-B_3(t-\tau)}\Pi_3 f(x(\tau))d\tau\\
&\quad-\int_{t}^\8e^{-B_1(t-\tau)}\Pi_1 f(x(\tau))d\tau.
\end{split}
\end{equation*}
We first verify that $\mathcal{T}$ maps $\mathscr{X}_{\b/2}$ into itself.

For notational convenience, we write $$0\wedge t=\min\{0,t\},\hs 0\vee t=\max\{0,t\}$$ for $t\in\R$.
Let $x\in \mathscr{X}_{\b/2}$. By \eqref{b1}-\eqref{b2} and \eqref{eC} we have
 \begin{equation}\label{e:2.15a}
\begin{split}
\|\mathcal {T}x(t)\|_\a&\leq Me^{\frac{\b}{4}|t|}\|y\|_\a+MC\int_{0\wedge t}^{0\vee t}e^{\frac{\b}{4}|t-\tau|}\big(\|x(\tau)\|_\a+1\big)d\tau\\
&\quad+ MC\int_{-\8}^t (t-\tau)^{-\a} e^{-\frac{3\b}{4}(t-\tau)}\big(\|x(\tau)\|_\a+1\big)d\tau\\
&\quad+MC\int_t^\8  e^{\frac{3\b}{4}(t-\tau)}\big(\|x(\tau)\|_\a+1\big)d\tau,
\end{split}\end{equation}where $C$ is the constant in \eqref{eC}.
Simple computations show  that
\begin{equation*}
\frac{\b}{4}|t-\tau|-\frac{\b}{2}|t|=-\frac{\b}{4}|t-\tau|-\frac{\b}{2}|\tau|,\Hs \tau\in [0\wedge t,\,0\vee t].
\end{equation*}
It is also easy to see  that
\begin{equation*}
\begin{split}
e^{-\frac{\b}{2}|t|}=e^{-\frac{\b}{2}|\tau+(t-\tau)|}\leq  e^{\frac{\b}{2}|t-\tau|}e^{-\frac{\b}{2}|\tau|},\Hs t,\,\tau\in\R.
\end{split}
\end{equation*}
Thus by \eqref{e:2.15a}  we find that
\begin{equation*}
\begin{split}
e^{-\frac{\b}{2}|t|}\|\mathcal {T}x(t)\|_\a
\leq\, &Me^{-\frac{\b}{4}|t|}\|y\|_\a+MC\int_{0\wedge t}^{0\vee t}e^{-\frac{\b}{4}|t-\tau|}\big[e^{-\frac{\b}{2}|\tau|}\big(\|x(\tau)\|_\a+1\big)\big]d\tau\\
&+ MC\int_{-\8}^t (t-\tau)^{-\a}e^{\frac{\b}{2}|t-\tau|} e^{-\frac{3\b}{4}(t-\tau)}\big[e^{-\frac{\b}{2}|\tau|}\big(\|x(\tau)\|_\a+1\big)\big]d\tau\\
&+MC\int_t^\8 e^{\frac{\b}{2}|t-\tau|}e^{\frac{3\b}{4}(t-\tau)}\big[e^{-\frac{\b}{2}|\tau|}\big(\|x(\tau)\|_\a+1\big)\big] d\tau\\
=\,&Me^{-\frac{\b}{4}|t|}\|y\|_\a+MC(\|x\|_{\mathscr{X}_{\b/2}}+1\big)\int_{0\wedge t}^{0\vee t}e^{-\frac{\b}{4}|t-\tau|}d\tau\\
&+MC(\|x\|_{\mathscr{X}_{\b/2}}+1\big)\int_{-\8}^t (t-\tau)^{-\a} e^{-\frac{\b}{4}(t-\tau)}d\tau\\
&+MC(\|x\|_{\mathscr{X}_{\b/2}}+1\big)\int_t^\8  e^{\frac{\b}{4}(t-\tau)} d\tau\\
\leq \,& M\|y\|_\a+M_\b C\big(\|x\|_{\mathscr{X}_{\b/2}}+1\big),\Hs \A\,t\in\R,
\end{split}
\end{equation*}
where \be\label{Mb} M_\b
 =M\int_0^\8\big(2+\tau^{-\a}\big)e^{-\frac{\b}{4}\tau}d\tau.\ee
It follows that  $\|\mathcal {T}x\|_{\mathscr{X}_{\b/2}}<\8$, that is, $\mathcal{T}x\in\mathscr{X}_{\b/2}$.

Next, we check that $\mathcal{T}$ is contractive.
Let $x,\,x'\in\mathscr{X}_{\b/2}$. In a quite similar fashion as above  it can be shown that
\begin{equation}\label{eq2.2}
\begin{split}
&\hs\,\, e^{-\frac{\b}{2}|t|}\|\mathcal {T}x(t)-\mathcal {T}x'(t)\|_\a\\
&\leq ML_f\int_{0\wedge t}^{0\vee t}e^{-\frac{\b}{4}|t-\tau|}\(e^{-\frac{\b}{2}|\tau|}\|x(\tau)-x'(\tau)\|_\a \)d\tau\\
&\quad+ ML_f\int_{-\8}^t(t-\tau)^{-\a} e^{-\frac{\b}{4}(t-\tau)}\(e^{-\frac{\b}{2}|\tau|}\|x(\tau)-x'(\tau)\|_\a \)d\tau\\
&\quad+ML_f\int^\8_t e^{\frac{\b}{4}(t-\tau)}\(e^{-\frac{\b}{2}|\tau|}\|x(\tau)-x'(\tau)\|_\a\)d\tau\\
&\leq M_\b L_f\|x-x'\|_{\mathscr{X}_{\b/2}},\Hs\A\, t\in\R,
\end{split}
\end{equation}
where (and below) $M_\b$ is the  number given in \eqref{Mb}.
Therefore
$$
\|\mathcal {T}x-\mathcal {T}x'\|_{\mathscr{X}_{\b/2}}\leq M_\b L_f\|x-x'\|_{\mathscr{X}_{\b/2}}.
$$
The condition \eqref{eq1.6} then asserts that $\cT$ is contractive.

Thanks to the Banach fixed-point theorem, $\cT=\cT_{\lam,y}$ has a unique fixed point $\gam_{y}:=\gam_{\lam,y}\in \mathscr{X}_{\b/2}$ which, by the definition of $\cT$, solves the integral equation
\begin{equation}\label{equ5}
\begin{split}
\gam_y(t)&=e^{-B_2t}y+\int_{0}^te^{-B_2(t-\tau)}\Pi_2 f(\gam_y(\tau))d\tau\\
&\quad+\int_{-\8}^{t}e^{-B_3(t-\tau)}\Pi_3 f(\gam_y(\tau))d\tau\\[0.5ex]
&\quad-\int_{t}^\8e^{-B_1(t-\tau)}\Pi_1 f(\gam_y(\tau))d\tau.
\end{split}
\end{equation}
(Hence $\gam_y(t)$   is  a solution of $(\ref{e1})$ on $\R$ with $\Pi_2\gam_y(0)=y$.)

Let $y,z\in X_2^\a$ and $t\in\R$. Similar to \eqref{eq2.2}, by  $(\ref{equ5})$ we find that
\begin{equation*}
\begin{split}
&e^{-\frac{\b}{2} |t|}\|\gam_{y}(t)-\gamma_{z}(t)\|_\a\\
\leq\,&Me^{-\frac{\b}{4} |t|}\|y-z\|_\a+ ML_f\int_{0\wedge t}^{0\vee t}e^{-\frac{\b}{4}|t-\tau|}\big(e^{-\frac{\b}{2}|\tau|}\|\gam_{y}(\tau)-\gamma_{z}(\tau)\|_\a \big)d\tau\\
&+ ML_f\int_{-\8}^t(t-\tau)^{-\a} e^{-\frac{\b}{4}(t-\tau)}\big(e^{-\frac{\b}{2}|\tau|}\|\gam_{y}(\tau)-\gamma_{z}(\tau)\|_\a \big)d\tau\\
&+ML_f\int^\8_t e^{\frac{\b}{4}(t-\tau)}\big(e^{-\frac{\b}{2}|\tau|}\|\gam_{y}(\tau)-\gamma_{z}(\tau)\|_\a \big)d\tau\\
\leq\,& M\|y-z\|_\a+ M_\b L_f\|\gam_{y}-\gamma_{z}\|_{\mathscr{X}_{\b/2}}.
\end{split}
\end{equation*}
Hence
$$
\|\gam_{y}-\gamma_{z}\|_{\mathscr{X}_{\b/2}}\leq\, M\|y-z\|_\a+ M_\b L_f\|\gam_{y}-\gamma_{z}\|_{\mathscr{X}_{\b/2}}.
$$
Therefore
 \begin{equation}\label{eq2}
\begin{split}
\|\gam_{y}(0)-\gamma_{z}(0)\|_\a\leq\|\gam_{y}-\gamma_{z}\|_{\mathscr{X}_{\b/2}}\leq \tilde L_0\|y-z\|_\a,\end{split}
\end{equation}
where $\tilde L_0={M}/({1-M_\b L_f})$ is a constant independent of $\lam\in J$.

Now we define a mapping $\xi_\lam: X_2^\a\ra X_{13}^\a$ as
$$
\xi_\lam(y)=\gam_y(0)-y,\Hs y\in X_2^\a.
$$
Setting $t=0$ in \eqref{equ5} one finds that
 \be\label{eq2.5}
\xi_\lam(y)=\int_{-\8}^0e^{B_3\tau}\Pi_3 f(\gamma_{y}(\tau))d\tau-\int_0^\8e^{B_1\tau}\Pi_1 f(\gamma_{y}(\tau))d\tau
\ee
for $y\in X_2^\a$. Let $L_0=\tilde L_0+1$. It follows by \eqref{eq2} that
  $$\|\xi_\lam(y)-\xi_\lam(z)\|_\a\leq  L_0\|y-z\|_\a, \Hs y,z\in X_2^\a.$$

 Let $$\cM^c_\lam=\{y+\xi_\lam(y):y\in X_2^\a\}.$$ Then $\cM^c_\lam$ is an invariant manifold of  \eqref{e:2.0}. Clearly $\cM_\lam^c$ is homeomorphic to $X_2^\a$. 
$\Box$

\Vs

\section{Bifurcation  and multiplicity of  \eqref{eq1}}
Let us now look at  the bifurcation and multiplicity of the Landesman-Lazer type problem of \eqref{eq1} when $\lam$ crosses any eigenvalue of the operator $A=-\Delta$.
For this purpose, we first turn our attention to the 
dynamic bifurcation of the  parabolic problem \be\label{e:3.1}\left\{\ba{ll}
  u_t-\Delta u=\lam u+f(x,u), \hs\,\, x\in
^{}  \Omega;\\[1ex]
  \,u(x)=0, \Hs\Hs \Hs \Hs x\in\partial\Omega,
\ea\right. \ee
where  $f\in C^1(\ol\Omega\X\R)$  is globally Lipschitz with Lipschitz constant $L_f$ and satisfies the (LC) condition in Section 1.
\subsection{Mathematical setting}

Let $H=L^2(\Omega)$ and $V=H_0^1(\Omega)$. By $(\cdot,\cdot)$ and
$|\cdot|$ we denote the usual inner product and norm on $H$,
respectively. The inner product and norm on $V$, denoted by $((\.,\.))$ and $\|\cdot\|$, respectively, are defined as
$$
((u,v))= \int_{\Omega}\nabla u\.\nab v\mathrm{d}x,\hs  \|u\|=\(\int_{\Omega}|\nabla u|^2\mathrm{d}x\)^{1/2}
$$
for $u,v\in V$. (The notation $\|\.\|$ is also used to denote the norm of any linear operator. We hope this will course no confusion.)

Denote $A$ the operator $-\Delta$ associated with the homogenous
Dirichlet boundary condition. Then $A$ is a sectorial operator on $H$ and
has a compact resolvent. It is a basic knowledge that $D(A^{1/2})=V$.

Let $$0<\mu_1<\mu_2<\cdots<\mu_k<\cdots$$ be the eigenvalues of $A$. Then $$\mb{$\mu_{k+1}-\mu_k\ra+\8$,\,\, as $k\ra\8$};$$ (see e.g. \cite[Chapter 4]{M}).
Denote $W_k$ the eigenspace corresponding to $\mu_k$.

System \eqref{e:3.1} can be rewritten as an abstract evolution equation  in $V$: \be\label{e3.2}
 u_t+B u= \tilde{f}(u),\hs u=u(t)\in V,\ee
 where $B:=B_\lam=A-\lam I,$ and $\tilde{f}:V\ra H$ is the Nemitski operator given by
 $$
 \tilde{f}(u)(x)=f(x,u(x)),\Hs u\in V.
 $$
One trivially verifies  that  \be\label{eq1.8}|\tilde{f}(u_1)- \tilde{f}(u_2)| \leq\tilde L \|u_1-u_2\|,\Hs\A\, u_1, u_2\in V,\ee
where $\tilde L =L_f/\sqrt{\mu_1}$.

\vs 
Set $$
\b_k=\min\{\mu_k-\mu_{k-1},\mu_{k+1}-\mu_k\},\Hs k=1,2,\cdots.
$$
(Here we assign $\b_1=\mu_2-\mu_1$.) For each $k$, if $\lam\in(\mu_k-\b_k,\mu_k+\b_k)$ then the  spectrum $\sig_\lam(B)$ of the operator $B=B_\lam$ has a spectrum decomposition $$\sig_\lam(B)=\sig_\lam^u\cup\sig_\lam^c\cup\sig_\lam^s,$$ where $\sig_\lam^c=\{\mu_k-\lam\} $, and $$\sig_\lam^u=\{\mu_j-\lam:\,\,j<k\},\hs
\sig_\lam^s=\{\mu_j-\lam:\,\,j>k\}.$$ Clearly
$$
\sig_\lam^u\subset (-\8,0),\hs\mb{and }\,\sig_\lam^s\subset (0,\8).
$$
The space $H$ has a corresponding orthogonal decomposition $H=\oplus_{i=u,c,s}H_i$ (independent of $\lam$) with $H_i\perp H_j$ if $i\ne j$.
 Set $H_{us}=H_u\oplus H_s$.

 Let $B_i$ be the restriction of $B$ on $H_i$, and denote $\Pi_i:H\ra H_i$   the projection, where $i=u,c,s,us$.


\subsection{Dynamic bifurcation from infinity}
 Denote $\Phi_\lam(t)$ the semiflow generated by the initial value problem of \eqref{e3.2} on $V$, namely, for each $u_0\in V$, $u(t)=\Phi_\lam(t)u_0$ is the (unique) strong  solution of  \eqref{e3.2} in $V$ with  $u(0)=u_0$.
Set $$V_i:=H_i\cap V,\Hs  i=u,s,c,us.$$ Denote $$M_{\b_k}=M\int_0^\8\big(2+\tau^{-\frac{1}{2}}\big)e^{-\frac{1}{4}\b_k\tau}d\tau,\Hs k\geq1.$$
Then we have the following result about the dynamic bifurcation  of  \eqref{e3.2}.

\bt\label{th3.1}Given $k\geq 1$, suppose  the Lipschitz constant $L_f$ of $f$ satisfies
 \be\label{e3b} M_{\b_k}L_f/\sqrt{\mu_1}<1.
\ee
   Then there exists $0<\theta <\b_k/4$ such that when $\lam\in \Lambda_k^-:=[\mu_k-\theta ,\mu_k)$, $\Phi_\lam$ has a compact invariant set $K^\8_\lam$ which takes the shape of $(m-1)$-dimensional sphere $\mathbb{S}^{m-1}$. Furthermore, \be\label{e3a}\lim_{\lam\ra\mu_k^-}\min\{\|v\|:\,\,v\in K^\8_\lam\}=\8.\ee\et

\br In the above theorem we have employed a topological concept, shape, without definition. Informally speaking, this notion can be seen as a generalization of that of homotopy type, and is used to describe topological structures of ``bad spaces'' such as invariant sets and attractors for which it is difficult to talk about homotopy type. It is a basic knowledge that spaces having  same homotopy type enjoy same shape.
The interested reader is referred to \cite{Kap} for details.
\er
\br As $\b_k\ra\8$ as $k\ra\8$, it is easy to see by definition that $M_{\b_k}\ra 0$ as $k\ra\8$. Therefore, for any globally Lipschitz function $f$, the smallness requirement \eqref{e3b} is automatically satisfied as long as  $k$ is large enough.
\er

To prove Theorem \ref{th3.1}, we need some  auxiliary results.
The proposition below is a straightforward application of Theorem \ref{th1}.
\bp\label{p2} Suppose $f$ satisfies \eqref{e3b}. Then for each $\lam\in(\mu_k-\b_k/4,\,\mu_k+\b_k/4):=J_k$,
the semiflow $\Phi_\lam$  has a global invariant   manifold $\cM^c_\lam$ given by \be\label{e3c}\cM^c_\lam=\{y+\xi_\lam(y):y\in V_c\},\ee where  $\xi_\lam:V_c\ra V_{13}$  is globally  Lipschitz with Lipschitz constant
\be\label{e3d}
L_{\xi_\lam}\leq L_0,\Hs\A\,\lam\in J_k
\ee
for some $L_0>0$ independent of $\lam\in J_k$.
\ep

Given a function $v$ on $\W$, we use $v_\pm$ to denote the positive and negtive parts of $v$, respectively, $$v_\pm=\max\{\pm v(x),0\},\hs x\in \W.$$ Then $v=v_+-v_-$.
The following  fundamental result on $f$ is taken from \cite{LLZ} (see also \cite[Section 6]{LiD}).
\bp\label{p1}\cite{LLZ} Suppose $f$ satisfies the Landesman-Lazer type condition $(LC)$.  Then for any $R,\ve>0$, there exists $s_0>0$ such that $$\int_\W f(x,sv+u)vdx\geq \int_\W \(\overline{f}v_++\underline{f}v_-\)dx-\ve$$
for all $s\geq s_0$, $v\in \overline{B} (1)$ and $u\in \overline{B} (R)$, where $B (r)$ denotes the ball in $H$ centered at $0$ with radius $r$.\ep

Henceforth we always assume  $L_f $ satisfies \eqref{e3b}, so  for each $\lam$ with $|\lam-\mu_k|<\b_k/4$, the semiflow $\Phi_\lam$ has an invariant  manifold $\cM_\lam^c$ given by \eqref{e3c}. If we reduce the system \eqref{e3.2}  on $\cM_\lam^c$, it takes the form
\be\label{eq2.1}w_t+B_c w=\Pi_2\tilde{f}(w+\xi_\lam(w)),\hs w\in  V_c.\ee
Let $\phi_\lam$ denote the semiflow generated by  \eqref{eq2.1} on $V_c$.

Since $V_c$ is finite dimensional, all the norms on $V_c$ are equivalent. Hence for convenience, we equip $V_c$ the norm $|\.|$ of $H$ in the following argument.

 Given $0\leq a<b\leq\8$, denote $$\Xi  [a,b]:=\{x\in V_c,\,\,a\leq|x| \leq b\}.$$
\bl\label{le1} Under the hypotheses \eqref{e3b} and $(LC)$, there exist $R_0\geq 0$ and $c_0>0$ such that the following assertions hold.
\begin{enumerate}
\item[$(1)$] If  $\lam\in[\mu_k,\mu_k+\b_k/4)$, then for any solution $w(t)$ of \eqref{eq2.1} in $\Xi [R_0,\8]$, we have
\be\label{e5.10}
\frac{d}{dt}|w| ^2\geq c_0 |w| .
\ee
\item[$(2)$] For any $R>R_0$, there exists $0<\ve<\b_k/4$ such that if  $\lam\in[\mu_k-\ve,\mu_k)$,
then \eqref{e5.10} holds true for any solution $w(t)$  of \eqref{eq2.1} in $\Xi  [R_0,R]$.
\item[$(3)$] There exists $\theta>0$ such that for each $\lam\in[\mu_k-\theta ,\mu_k)$, the semiflow $\phi_\lam$ has a positively invariant  set $N_\lam=\Xi [a_\lam,b_\lam]$ with $a_\lam,\,b_\lam\ra\8$ as $\lam\ra \mu_k^-$. \end{enumerate}
 \el
\vs
{\bf Proof.}
Taking the inner product of (\ref{eq2.1}) with $w$ in $H$, it yields
\be\label{e6.16}\frac{1}{2}\frac{d}{dt}|w| ^2+(\mu_k-\lam)|w| ^2= \big(\Pi_2\tilde{f}(w+\xi_\lam(w)),w\big)= \big(\tilde{f}(w+\xi_\lam(w)),w\big).
\ee
Let us first estimate the last term  in
\eqref{e6.16}.

As the  norm $|\.|_{L^1(\W)}$ of $L^1(\W)$ and that of  $H$ are equivalent on $V_c$, one easily sees that
  \be\label{e6.26}\ba{ll}
\min\{|v|_{L^1(\Omega)}:\,\,v\in V_c,\,\,|v| =1\}:=r>0.\ea
\ee
Pick a number $\de>0$ with $\de\leq\min\{\ol f,\ul f\}$. We infer from the representation of $\xi_\lam$ and the boundedness of $f$ that $\xi_\lam$ is uniformly bounded on $\lam$.
Thus by virtue of Proposition \ref{p1} there exists $s_0>0$  such that when $s\geq s_0$,
\be\label{e6.20}\begin{split}\big(\tilde{f}(sv+\xi_\lam(sv)),v\big)&=\int_\Omega f(x,sv+\xi_\lam(sv))v\,dx\\[1ex]
&\geq \int_\Omega\(\ol fv_++\ul fv_-\)dx-\frac{1}{2}r\de \end{split}\ee for all $v\in \ol\mB (1)$.

Now we rewrite $$\mb{ $w=sv$,\,\, where
$s=|w| $.}$$ Then $v\in \ol\mB (1)$. Suppose $s\geq s_0$, by
(\ref{e6.20}) one finds that
\begin{equation*}\begin{split}
\big(\tilde{f}\big(w+\xi_\lam(w)),w\big)&=s\int_\W f\(x,sv+\xi_\lam(sv)\)vdx\\[1ex]
&\geq s\(\int_\Omega\(\ol f v_++\ul
f v_-\)dx-\frac{1}{2}r\de \).\end{split}
\end{equation*}
 Observing that
\begin{equation*}\begin{split}
&\int_\Omega\(\ol fv^++\ul fv^-\)dx -\frac{1}{2}r\de \\[1ex]
\geq &\de \int_\Omega|v|dx-\frac{1}{2}r\de\geq (\mb{by }(\ref{e6.26}))\geq \frac{1}{2}r\de,\end{split}\end{equation*}
we conclude that \be\label{e5.22} \big(\tilde{f}\big(w+\xi_\lam(w)),w\big)\geq
\frac{1}{2}r\de s=\frac{1}{2}r\de |w| . \ee

Combining \eqref{e5.22} and (\ref{e6.16}) together we find that
\be\label{e5.23}
\frac{d}{dt}|w| ^2\geq 2\(\lam-\mu_k\) |w| ^2+r\de |w|
\ee as long as $|w(t)| \geq s_0$.

\vs
Set  $R_0=s_0$,  $c_0=r\de/2$. Assume $\lam\in[\mu_k,\mu_k+\b_k/4)$. Then $\lam-\mu_k\geq 0$, and we infer from \eqref{e5.23} that if $|w| \geq R_0$ then
$$
\frac{d}{dt}|w| ^2\geq r\de |w| >c_0|w| .$$
 Hence assertion  (1) holds.

\vs
Now assume $\lam<\mu_k$. Let $R>R_0$. Choose an $\ve>0$ sufficiently small so that   $\ve R^2<r\de s_0/4$. Then if $\lam\in[\mu_k-\ve,\mu_k)$, for any solution $w(t)$  of \eqref{eq2.1} in $\Xi  [R_0,R]$, by \eqref{e5.23} we deduce that
$$\ba{ll}
\frac{d}{dt}|w| ^2&\geq -2|\lam-\mu_k|\,R^2+r\de |w| \\[1ex]
&\geq -2\theta R^2 +2c_0|w| \\[1ex]
&= \big(c_0|w| -2\theta R^2\big)+c_0|w| \\
&\geq c_0|w| ,\ea
$$  which justifies  assertion (2).

Note that \eqref{e5.10} implies that  $\Xi [R,\8]$ is positively invariant for $\phi_\lam$ when $\lam\in[\mu_k-\ve,\mu_k)$.
\vs
 Let $R_j=R_0+j$ ($j=1,2,\cdots$). Then for each $j$, we deduce by assertion (2) that there exists $\ve_j>0$ such that if $\lam\in[\mu_k-\ve_j,\mu_k)$, \eqref{e5.10} holds true for any solution $w(t)$  of \eqref{eq2.1} in $\Xi  [R_0,R_j]$. Hence $\Xi [R_j,\8]$ is positively invariant for $\phi_\lam$ when $\lam\in[\mu_k-\ve_j,\mu_k)$.

On the other hand, by the boundedness of $f$ we have
\begin{equation}\label{eq3.6}
\begin{split}\big(\tilde{f}(w+\xi_\lam(w)),w\big)&\leq \big|\tilde{f}(w+\xi_\lam(w))\big| \,|w| \leq C\,|w| \\
&\leq \frac{\mu_k-\lam}{2}|w|^2 +C(\lam),\end{split}
\end{equation}
where  $C(\lam)\ra +\8$ as $\lam\ra \mu_k^-$.
Combining \eqref{eq3.6} with \eqref{e6.16}, we find that
$$\frac{d}{dt}|w| ^2\leq-(\mu_k-\lam)|w| ^2+2C(\lam).
$$
Thanks to the classical Gronwall lemma,
\be\label{e3e}|w(t)| ^2\leq e^{-(\mu_k-\lam)t}|w(0)| ^2+\(1-e^{-(\mu_k-\lam)t}\)\frac{2C(\lam)}{\mu_k-\lam},\hs t\geq0.\ee
By \eqref{e3e} it is easy to verify that if $$\rho\geq \sqrt{2C(\lam)/(\mu_k-\lam)}:=\rho_\lam,$$ then $\{v\in V_c:\,\,|v|\leq \rho\}$ is positively invariant.

We may assume  $\ve_1>\ve_2>\cdots>\ve_j\ra 0$. Then $$[\mu_k-\ve_1,\mu_k)=\Cup_{j\geq 1}[\mu_k-\ve_j,\,\mu_k-\ve_{j+1}).$$ Set $\theta =\ve_1$, and let $\lam\in[\mu_k-\theta ,\mu_k)$. If $\lam\in [\mu_k-\ve_j,\,\mu_k-\ve_{j+1})$,  we take
$$
a_\lam=R_j,\hs b_\lam=R_j+\rho_\lam.
$$
Clearly $a_\lam,b_\lam\ra\8$ as $\lam\ra\mu_k^-$. We infer from the above argument that $\Xi[a_\lam,b_\lam]$ is  positively invariant under the system $\phi_\lam$, hence assertion (3) holds true. $\Box$
\Vs
Now let us turn  to the proof of  Theorem \ref{th3.1}.
\vs\noindent {\bf Proof of Theorem \ref{th3.1}.}
Let $k\geq1$, and let $\theta $ be the number given in Lemma \ref{le1}. Assume $\lam\in[\mu_k-\theta ,\mu_k)$. Then  Lemma \ref{le1} (3) asserts that $N_\lam=\Xi[a_\lam,b_\lam]$ is a positively invariant set of   $\phi_\lam$.  Set
$$
\cA_\lam^\8=\Cap_{\tau\geq 0}\ol{\Cup_{t\geq\tau}\phi_\lam(t)N_\lam}.
$$
By the basic knowledge in the attractor theory (see e.g. \cite{Hale,Tem}) we know that $\cA_\lam^\8$ is the global attractor of $\phi_\lam$ restricted on the phase space $X=N_\lam$.

Since $N_\lam$ has the homotopy type of an  $(m-1)$-dimensional sphere $\mathbb{S}^{m-1}$, it shares the same shape of $\mathbb{S}^{m-1}$. Therefore employing the shape theory of attractors in \cite{Kap} (see also \cite{S}), $\cA^\8_\lam$ has the shape of $\mathbb{S}^{m-1}$.

  Let $$K^\8_\lam=\{w+\xi_\lam(w): w\in \cA_\lam^\8\}. $$
Then  $K^\8_\lam\subset \cM^c_\lam$ is a compact invariant set of the original system $\Phi_\lam$ which takes the shape of an $(m-1)$-dimensional sphere.

The conclusion in \eqref{e3a} directly follows from the fact that  $a_\lam,b_\lam\ra+\8$ as $\lam\ra \mu_k^-$. The proof of the theorem is complete.  $\Box$
\subsection{Bifurcation and multiplicity of \eqref{eq1}}
We are now ready  to state and prove our main result in this work.

\bt\label{th3.2} Assume $f$ satisfies the Landesman-Lazer type condition $(LC)$.
Let $\mu_k$ be an eigenvalue of $A$. Suppose  $L_f$ satisfies  \eqref{e3b}. Then there exists $0<\theta < \b_k/4$ such that the problem \eqref{eq1} has at least three distinct solutions  $u_{\lam}^{i}$ $(i=0,1,2)$ for each $\lam\in\Lam_k^-=[\mu_k-\theta ,\mu_k)$, and $u_{\lam}^{1}$ and $u_{\lam}^{2}$ go to $\8$ as $\lam\ra \mu_k^-$ whereas $u_\lam^0$ remains bounded on $\Lam^-_k$.
 \et
 \br It is known under the hypotheses of Theorem \ref{th3.2} that \eqref{eq1} has at least one solution for all $\lam\in\R$; see, e.g.,
 \cite[Theorem 5.3]{LLZ}.
 \er
 {\bf Proof.} Let  $\mu_k$ be an eigenvalue of $A$,  and let $\theta $ be the number given in Theorem \ref{th3.1}. Then Theorem \ref{th3.1} asserts that for each  $\lam\in \Lambda_k^-=[\mu_k-\theta ,\mu_k)$, the system $\Phi_\lam$ has a compact invariant set  $K^\8_\lam$. We first show that $K^\8_\lam$ contains at least two distinct solutions
  $u_{\lam}^{1}$ and $u_{\lam}^{2}$ of \eqref{eq1}.

  Since $K^\8_\lam$ has the shape of $\mathbb{S}^{m-1}$, it consists of at least two distinct points $u$ and $v$. If both $u$ and $v$ are equilibrium points of $\Phi_\lam$ then we are done with $u_{\lam}^{1}=u$ and $u_{\lam}^{2}=v$. Thus we assume, say, $v$ is not
 an equilibrium of $\Phi_\lam$. Let $\gam=\gam(t)$ be a complete solution of $\Phi_\lam$  contained in $K^\8_\lam$, $\gam(0)=v$. As $\Phi_\lam$ is a gradient system (and hence there is no homoclinic structure in $K^\8_\lam$),  the $\omega $-limit $\omega (\gam)$ and $\a$-limit $\a(\gam)$ of $\gam$ do not intersect. Because $\omega (\gam)$ and $\a(\gam)$ consist of equilibrium points, we deduce that $\Phi_\lam$ has at least two distinct equilibria in $K^\8_\lam$. Consequently \eqref{eq1} has two distinct solutions $u_{\lam}^{1}$ and $u_{\lam}^{2}$.

By virtue of \eqref{e3a} it is clear that
\be\label{eq3.7}\lim_{\lam\ra\mu_k^-}\|u_{\lam}^{i}\|=\8,\Hs i=1,2.
\ee

We also infer from  \cite[Theorem 5.3]{LLZ} that for each $\lam\in\Lam^-_k$, $\Phi_\lam$ has an equilibrium $u^0_\lam$ which remains bounded as $\lam$ varies in $\Lam^-_k$. In consideration of \eqref{eq3.7} it can be assumed that $u^0_\lam\ne u^i_\lam$ for $i=1,2$. Then $u^i_\lam$ ($i=0,1,2$) are solutions of \eqref{eq1} fulfilling the requirements of the theorem. $\Box$

\Vs \noindent{\bf Acknowledgments.} \vs
The authors gratefully acknowledge the support of the National Natural Science Foundation of China under Grants 11471240, 11071185.

{\small

\end{document}